\documentclass{amsart}
\usepackage{amsmath}
\usepackage{amssymb}
\usepackage{amsfonts}

\newtheorem{theorem}{Theorem}
\theoremstyle{plain}

\newtheorem{definition}{Definition}

\newtheorem{lemma}{Lemma}

\newtheorem{remark}{Remark}

\numberwithin{equation}{section}
\input{tcilatex}

\begin{document}
\title[On new weighted inequalities]{On weighted Montogomery identities for
Riemann-Liouville fractional integrals }
\author{Mehmet Zeki Sarikaya$^{\star }$}
\address{Department of Mathematics, Faculty of Science and Arts, D\"{u}zce
University, D\"{u}zce, Turkey}
\email{sarikayamz@gmail.com}
\thanks{$^{\star }$corresponding author}
\author{Hatice YALDIZ}
\email{yaldizhatice@gmail.com}

\begin{abstract}
In this paper, we extend the weighted Montogomery identities for the
Riemann-Liouville fractional integrals. We also use this Montogomery
identities to establish some new Ostrowski type integral inequalities.
\end{abstract}

\maketitle

\section{Introduction}

The inequality of Ostrowski \cite{Ostrowski} gives us an estimate for the
deviation of the values of a smooth function from its mean value. More
precisely, if $f:[a,b]\rightarrow \mathbb{R}$ is a differentiable function
with bounded derivative, then%
\begin{equation*}
\left\vert f(x)-\frac{1}{b-a}\int\limits_{a}^{b}f(t)dt\right\vert \leq \left[
\frac{1}{4}+\frac{(x-\frac{a+b}{2})^{2}}{(b-a)^{2}}\right] (b-a)\left\Vert
f^{\prime }\right\Vert _{\infty }
\end{equation*}
for every $x\in \lbrack a,b]$. Moreover the constant $1/4$ is the best
possible.

For some generalizations of this classic fact see the book \cite[p.468-484]%
{mitrovich} by Mitrinovic, Pecaric and Fink. A simple proof of this fact can
be \ done by using the following identity \cite{mitrovich}:

If $f:[a,b]\rightarrow \mathbb{R}$ is differentiable on $[a,b]$ with the
first derivative $f^{\prime }$ integrable on $[a,b],$ then Montgomery
identity holds:%
\begin{equation}
f(x)=\frac{1}{b-a}\int\limits_{a}^{b}f(t)dt+\int\limits_{a}^{b}P_{1}(x,t)f^{%
\prime }(t)dt,  \label{h}
\end{equation}%
where $P_{1}(x,t)$ is the Peano kernel defined by%
\begin{equation*}
P_{1}(x,t):=\left\{ 
\begin{array}{ll}
\dfrac{t-a}{b-a}, & a\leq t<x \\ 
&  \\ 
\dfrac{t-b}{b-a}, & x\leq t\leq b.%
\end{array}%
\right.
\end{equation*}%
Recently, several generalizations of the Ostrowski integral inequality are
considered by many authors; for instance covering the following concepts:
functions of bounded variation, Lipschitzian, monotonic, absolutely
continuous and $n$-times differentiable mappings with error estimates with
some special means together with some numerical quadrature rules. For recent
results and generalizations concerning Ostrowski's inequality, we refer the
reader to the recent papers \cite{Cerone1}, \cite{Duo}, \cite{Dragomir}-\cite%
{Liu}, \cite{sarikaya}-\cite{sarikaya2}.

In \cite{Anastassiou} and \cite{sarikaya3}, the authors established some
inequalities for differentiable mappings which are connected with Ostrowski
type inequality by used the Riemann-Liouville fractional integrals, and they
used the following lemma to prove their results:

\begin{lemma}
Let $f:I\subset \mathbb{R}\rightarrow \mathbb{R}$ be differentiable function
on $I^{\circ }$ with $a,b\in I$ ($a<b$) and $f^{\prime }\in L_{1}[a,b]$, then%
\begin{equation}
f(x)=\frac{\Gamma (\alpha )}{b-a}(b-x)^{1-\alpha }{\Large J}_{a}^{\alpha
}f(b)-{\Large J}_{a}^{\alpha -1}(P_{2}(x,b)f(b))+{\Large J}_{a}^{\alpha
}(P_{2}(x,b)f^{^{\prime }}(b)),\ \ \ \alpha \geq 1,  \label{z}
\end{equation}%
where $P_{2}(x,t)$ is the fractional Peano kernel defined by%
\begin{equation*}
P_{2}(x,t)=\left\{ 
\begin{array}{ll}
\dfrac{t-a}{b-a}(b-x)^{1-\alpha }\Gamma (\alpha ), & a\leq t<x \\ 
&  \\ 
\dfrac{t-b}{b-a}(b-x)^{1-\alpha }\Gamma (\alpha ), & x\leq t\leq b.%
\end{array}%
\right.
\end{equation*}
\end{lemma}

In this article, we use the Riemann-Liouville fractional integrals to
establish some new weighted integral inequalities of Ostrowski's type. From
our results, the weighted and the classical Ostrowski's inequalities can be
deduced as some special cases.

\section{Fractional Calculus}

Firstly, we give some necessary definitions and mathematical preliminaries
of fractional calculus theory which are used further in this paper. More
details, one can consult \cite{gorenflo}, \cite{samko}.

\begin{definition}
The Riemann-Liouville fractional integral operator of order $\alpha \geq 0$
with $a\geq 0$ is defined as%
\begin{eqnarray*}
J_{a}^{\alpha }f(x) &=&\frac{1}{\Gamma (\alpha )}\dint\limits_{a}^{x}(x-t)^{%
\alpha -1}f(t)dt, \\
J_{a}^{0}f(x) &=&f(x).
\end{eqnarray*}
\end{definition}

Recently, many authors have studied a number of inequalities by used the
Riemann-Liouville fractional integrals, see (\cite{Anastassiou}, \cite%
{Belarbi}, \cite{Dahmani}, \cite{Dahmani1}, \cite{sarikaya3}, \cite%
{sarikaya4}) and the references cited therein.

\section{Main Results}

Throughout this work, we assume that the weight function $w:\left[ a,b\right]
\rightarrow \lbrack 0,\infty ),$ is integrable, nonnegative and 
\begin{equation*}
m(a,b)=\int\limits_{a}^{b}w(t)dt<\infty .
\end{equation*}

In order to prove our main results, we need the following identities:

\begin{lemma}
\label{lm} Let $f:I\subset \mathbb{R}\rightarrow \mathbb{R}$ be a
differentiable function on $I^{\circ }$ with $a,b\in I$ ($a<b$)$,$ $\alpha
\geq 1$ and $f^{\prime }\in L_{1}[a,b]$, then the generalization of the
weighted Montgomery identity for fractional integrals holds:
\end{lemma}

\begin{eqnarray}
m\left( a,b\right) f\left( x\right) &=&\left( b-x\right) ^{1-\alpha }\Gamma
\left( \alpha \right) J_{a}^{\alpha }\left( w\left( b\right) f\left(
b\right) \right)  \notag \\
&&  \label{1} \\
&&-J_{a}^{\alpha -1}\left( \Omega _{w}\left( x,b\right) f\left( b\right)
\right) +J_{a}^{\alpha }\left( \Omega _{w}\left( x,b\right) f^{^{\prime
}}\left( b\right) \right)  \notag
\end{eqnarray}%
where $\Omega _{w}\left( x,t\right) $ is the weighted fractional Peano
kernel defined by

\begin{equation}
\Omega _{w}\left( x,t\right) :=\left\{ 
\begin{array}{ll}
\left( b-x\right) ^{1-\alpha }\Gamma \left( \alpha \right)
\dint\limits_{a}^{t}w\left( u\right) du, & t\in \lbrack a,x) \\ 
\left( b-x\right) ^{1-\alpha }\Gamma \left( \alpha \right)
\dint\limits_{b}^{t}w\left( u\right) du, & t\in \lbrack x,b].%
\end{array}%
\right.  \label{2}
\end{equation}

\begin{proof}
By definition of $\Omega _{w}\left( x,t\right) $, we have%
\begin{eqnarray}
&&J_{a}^{\alpha }\left( \Omega _{w}\left( x,b\right) f^{^{\prime }}\left(
b\right) \right)  \notag \\
&&  \label{3} \\
&=&\frac{1}{\Gamma \left( \alpha \right) }\dint\limits_{a}^{b}\left(
b-t\right) ^{\alpha -1}\Omega _{w}\left( x,t\right) f^{\prime }\left(
t\right) dt  \notag \\
&&  \notag \\
&=&\left( b-x\right) ^{1-\alpha }\left[ \dint\limits_{a}^{x}\left(
b-t\right) ^{\alpha -1}\left( \dint\limits_{a}^{t}w\left( u\right) du\right)
f^{^{\prime }}\left( t\right) dt\right.  \notag \\
&&  \notag \\
&=&\left. +\dint\limits_{x}^{b}\left( b-t\right) ^{\alpha -1}\left(
\dint\limits_{b}^{t}w\left( u\right) du\right) f^{^{\prime }}\left( t\right)
dt\right]  \notag \\
&&  \notag \\
&=&\left( b-x\right) ^{1-\alpha }\left( J_{1}+J_{2}\right) .  \notag
\end{eqnarray}%
Integrating by parts, we can state:%
\begin{eqnarray}
J_{1} &=&\left( b-x\right) ^{\alpha -1}\left( \dint\limits_{a}^{x}w\left(
u\right) du\right) f\left( x\right)  \notag \\
&&  \label{41} \\
&&+\left( \alpha -1\right) \dint\limits_{a}^{x}\left( b-t\right) ^{\alpha
-2}\left( \dint\limits_{a}^{t}w\left( u\right) du\right) f\left( t\right)
dt-\dint\limits_{a}^{x}\left( b-t\right) ^{\alpha -1}w\left( t\right)
f\left( t\right) dt  \notag
\end{eqnarray}%
and similary,%
\begin{eqnarray}
J_{2} &=&\left( b-x\right) ^{\alpha -1}\left( \dint\limits_{x}^{b}w\left(
u\right) du\right) f\left( x\right)  \notag \\
&&  \label{5} \\
&&+\left( \alpha -1\right) \dint\limits_{x}^{b}\left( b-t\right) ^{\alpha
-2}\left( \dint\limits_{b}^{t}w\left( u\right) du\right) f\left( t\right)
dt-\dint\limits_{x}^{b}\left( b-t\right) ^{\alpha -1}w\left( t\right)
f\left( t\right) dt.  \notag
\end{eqnarray}%
Adding (\ref{41}) and (\ref{5}), we obtain (\ref{1}) which this completes
the proof.
\end{proof}

\begin{remark}
If we choose $\alpha =1$ and $w\left( u\right) =1$, the formula (\ref{1})
reduces to the classical Montgomery Identity given by (\ref{h}).
\end{remark}

\begin{remark}
If we choose $w\left( u\right) =1$, the formula (\ref{1}) reduces to the
fractional Montgomery Identity given by (\ref{z}).
\end{remark}

\begin{theorem}
\label{thm3} Let $f:[a,b]\rightarrow \mathbb{R}$ be differentiable on $(a,b)$
such that $f^{^{\prime }}\in L_{1}[a,b],$ where $a<b.$ If $\left\vert
f^{^{\prime }}(x)\right\vert \leq M$ for every $x\in \lbrack a,b]$ and $%
\alpha \geq 1$, then the following Ostrowski fractional inequality holds:
\end{theorem}

\begin{eqnarray}
&&\left\vert m\left( a,b\right) f\left( x\right) -\left( b-x\right)
^{1-\alpha }\Gamma \left( \alpha \right) J_{a}^{\alpha }\left( w\left(
b\right) f\left( b\right) \right) -J_{a}^{\alpha -1}\left( \Omega _{w}\left(
x,b\right) f\left( b\right) \right) \right\vert  \notag \\
&&  \label{22} \\
&\leq &\frac{M\left( b-x\right) ^{1-\alpha }}{\alpha }\left[ A(x)-\left(
b-x\right) ^{\alpha }B(x)\right]  \notag
\end{eqnarray}%
where%
\begin{equation*}
A(x)=\dint\limits_{a}^{x}\left( b-u\right) ^{\alpha -1}w\left( u\right)
du-\dint\limits_{x}^{b}\left( b-u\right) ^{\alpha }w\left( u\right) du
\end{equation*}%
and 
\begin{equation*}
B(x)=\dint\limits_{a}^{x}w\left( u\right) du-\dint\limits_{x}^{b}w\left(
u\right) du.
\end{equation*}

\begin{proof}
From Lemma \ref{lm}, we get%
\begin{eqnarray}
&&\left\vert m\left( a,b\right) f\left( x\right) -\Gamma \left( \alpha
\right) \left( b-x\right) ^{1-\alpha }J_{a}^{\alpha }\left( w\left( b\right)
f\left( b\right) \right) -J_{a}^{\alpha -1}\left( \Omega _{w}\left(
x,b\right) f\left( b\right) \right) \right\vert  \notag \\
&&  \notag \\
&\leq &\frac{1}{\Gamma \left( \alpha \right) }\left\vert
\dint\limits_{a}^{b}\left( b-t\right) ^{\alpha -1}\Omega _{w}\left(
x,t\right) f^{^{\prime }}\left( t\right) dt\right\vert  \notag \\
&&  \label{23} \\
&\leq &\frac{M}{\Gamma \left( \alpha \right) }\dint\limits_{a}^{b}\left(
b-t\right) ^{\alpha -1}\left\vert \Omega _{w}\left( x,t\right) \right\vert dt
\notag \\
&&  \notag \\
&=&M\left( b-x\right) ^{1-\alpha }\left( \dint\limits_{a}^{x}\left(
b-t\right) ^{\alpha -1}\left( \dint\limits_{a}^{t}w\left( u\right) du\right)
dt+\dint\limits_{x}^{b}\left( b-t\right) ^{\alpha -1}\left(
\dint\limits_{t}^{b}w\left( u\right) du\right) dt\right)  \notag \\
&&  \notag \\
&=&M(b-x)^{1-\alpha }\left\{ J_{3}+J_{4}\right\} .  \notag
\end{eqnarray}%
Now, using the change of order of integration we get%
\begin{eqnarray*}
J_{3} &=&\dint\limits_{a}^{x}\left( b-t\right) ^{\alpha -1}\left(
\dint\limits_{a}^{t}w\left( u\right) du\right) dt \\
&& \\
&=&\dint\limits_{a}^{x}w(u)\dint\limits_{u}^{x}\left( b-t\right) ^{\alpha
-1}dtdu \\
&& \\
&=&\frac{1}{\alpha }\left[ \dint\limits_{a}^{x}\left( b-u\right) ^{\alpha
-1}w\left( u\right) du-\left( b-x\right) ^{\alpha
}\dint\limits_{a}^{x}w\left( u\right) du\right]
\end{eqnarray*}%
and similarly,%
\begin{eqnarray*}
J_{4} &=&\dint\limits_{x}^{b}\left( b-t\right) ^{\alpha -1}\left(
\dint\limits_{t}^{b}w\left( u\right) du\right) dt \\
&& \\
&=&\dint\limits_{x}^{b}w(u)\dint\limits_{x}^{u}\left( b-t\right) ^{\alpha
-1}dtdu \\
&& \\
&=&\frac{1}{\alpha }\left[ \left( b-x\right) ^{\alpha
}\dint\limits_{x}^{b}w\left( u\right) du-\dint\limits_{x}^{b}\left(
b-u\right) ^{\alpha }w\left( u\right) du\right] .
\end{eqnarray*}%
Using $J_{3}$ and $J_{4}$ in (\ref{23}), we obtain (\ref{22}).
\end{proof}

\begin{remark}
We note that in the special cases, if we take $w\left( u\right) =1$ in
Theorem \ref{thm3}, then it reduces Theorem 4.1 proved by Anastassiou et.
al. \cite{Anastassiou}. So, our results are generalizations of the
corresponding results of Anastassiou et. al. \cite{Anastassiou}.
\end{remark}

\end{document}